\documentclass[11pt]{article}
\usepackage{amsmath,amssymb,amsfonts}
\usepackage[dvips]{graphicx}
\usepackage{mathrsfs}

\newtheorem{theorem}{Theorem}[section]

\newcommand{\qed}{\hfill \ensuremath{\Box}}
\DeclareMathOperator{\cosec}{cosec}

\numberwithin{equation}{section}

\begin{document}

\title{\Large      More Jordan type inequalities 		}
\author{        D. Aharonov and U. Elias }
\date{}  
\maketitle


\begin{abstract}
The function  $ \tan( \pi x / 2) / ( \pi x / 2 ) $  is expanded into a Laurent series of  
$ 1 - x^2 $,  where the coefficients are given explicitly as combinations of zeta function 
of even integers.
This is used to achieve a sequence of upper and lower bounds which are very precise even 
at the poles at $ x = \pm 1 $.

Similar results are obtained for other trigonometric functions with poles.

\end{abstract}

\vspace{-3 mm}


\section{ Introduction }
\setcounter{section}{1}	  

Jordan's inequality,
$$
	\frac{2}{\pi}	\le    \frac {\sin x} {x}  \le	1,	
		\qquad		-\frac{\pi}{2} \le x \le \frac{\pi}{2}, 
$$
has numerous extensions and generalizations.  (To simplify the notation we agree that 
 for $ x=0 $ the value of  $ \sin x / x $  is  1).   An exhaustive review of the literature
is available in  \cite{Qi-review}.  One of the many types of generalizations is a sequence 
of the two sided inequalities.  The first of them is the inequality 
\begin{equation}					
	\frac{2}{\pi} + \frac{1}{\pi^3} (\pi^2 - 4 x^2 )		\label{eq:J-1} 
\le	\frac {\sin x} {x}  
\le	\frac{2}{\pi} + \frac{\pi - 2}{\pi^3} (\pi^2 - 4 x^2 ) , 
	\qquad		-\frac{\pi}{2} \le x \le \frac{\pi}{2}   
\end{equation}
by Ling Zhu  \cite{Zhu1}. 
It is followed by 
\begin{equation}							\label{eq:J-2} 
	\frac{2}{\pi} + \frac{1}{\pi^3} (\pi^2 - 4 x^2 ) + \frac{12- \pi^2}{16\pi^5} (\pi^2 - 4 x^2 )^2
\le	\frac {\sin x} {x}  
\le	\frac{2}{\pi} + \frac{\pi - 2}{\pi^3} (\pi^2 - 4 x^2 ) + \frac{\pi-3}{\pi^5} (\pi^2 - 4 x^2 )^2,	
\end{equation}
by  \cite{Zhu2}; 
Finally Ling Zhu verifies  \cite{Zhu3}  the general 
\begin{equation}							\label{eq:J-N} 
	\sum_{n=0}^N a_n (\pi^2 - 4 x^2 )^n + a_{N+1} (\pi^2 - 4 x^2 )^{N+1}
\le	\frac {\sin x} {x}  
\le	\sum_{n=0}^N a_n (\pi^2 - 4 x^2 )^n + b_{N+1} (\pi^2 - 4 x^2 )^{N+1},	
\end{equation}
where $ a_{n} $  are given by a recursion formula.   See also \cite[Section 3.7]{Qi-review}.  

Later  (\ref{eq:J-N})  was extended in  \cite{Zhu4}  to 
\begin{equation}							\label{eq:J-r}  
	\sum_{n=0}^N a_{2n,r} (r^2 - x^2 )^n + \alpha_{N+1,r} (r^2 - x^2 )^{N+1}
\le	\frac {\sin x} {x}  
\le	\sum_{n=0}^N a_{2n,r} (r^2 - x^2 )^n + \beta_{N+1,r} (r^2 - x^2 )^{N+1} 
\end{equation} 
for  $ |x| \le r \le \pi / 2 $.  The proofs of inequalities  (\ref{eq:J-1})--(\ref{eq:J-r})  
are based on a monotone version of the l'Hospital rule:  
	{\it  If  $ g'(x) \neq 0 $  and  $ f'(x) / g'(x) $  is increasing on  $ (a,b) $, so is the 
	quotient  $ ( f(x) - f(a) ) / ( g(x) - g(a) ) $}.

The two sided inequalities  (\ref{eq:J-N})  are closely related to the infinite expansion 
of   \cite{Jian-Lin Li} 
\begin{equation}							\label{eq:J-series}
\frac {\sin x} {x}  
	= \frac{2}{\pi} + \sum_{n=1}^\infty \frac{ (-1)^n R_n } { n!\pi^{2n} }  (\pi^2 - 4 x^2 )^n  
\end{equation} 
where 
$ 
\displaystyle   
R_k 
= \sum_{ n=k }^\infty \frac{ (-1)^n n! } { (2n+1)! (n-k)! } \left( \frac{ \pi } {2} \right)^{2n} 
$
and  $ (-1)^k R_k > 0 $.   See also  \cite[Section 2.3]{Qi-review}.
Identity  (\ref{eq:J-series})  is verified by a straightforward rearrangement of a double 
power series.  By direct integration one achieves also 
\begin{equation*}
\cos x 
	= \frac{1}{4\pi} (\pi^2 - 4 x^2 ) 
	  + \sum_{n=1}^\infty \frac{ (-1)^n R_n } { 8(n+1)!\pi^{2n} }  (\pi^2 - 4 x^2 )^{n  +1} .
\end{equation*}

Similar ideas were developed for Bessel functions in a paper by Baricz and Wu  \cite{Baricz-Wu}
and are summarized in  \cite{Baricz}.  In  \cite[p. 145]{Baricz}  the author uses the notation  
	$ \mathscr{J}_p(x) = 2^p \Gamma(p+1) x^{-p} J_p(x) $,  
normalized so that  $ \mathscr{J}_p(0) = 1 $.
Without normalization the result of Baricz and Wu are formulated as 
\begin{equation}							\label{eq:Baricz} 
\begin{aligned}
\sum_{k=0}^N  \frac{ r^{-(p+k) } J_{p+k}(r) } { 2^k k! }  & (r^2 - x^2 )^k 
		+ \alpha_{p} (r^2 - x^2 )^{N+1}  \  \ 
    \le	 \   x^{-p} J_p(x)			\\
&   \le	 \ 
\sum_{k=0}^N  \frac{ r^{-(p+k) } J_{p+k}(r) } { 2^k k! } (r^2 - x^2 )^k 
		+ \beta_{p} (r^2 - x^2 )^{N+1} 
\end{aligned}
\end{equation}
for all  $ N $  and for  $ |x| \le r $,  where  $ 0 < r \le j_{p+1,1} $  and  $ j_{p+1,1} $  
denotes the first positive zero of  $ J_{p+1}(x) $.  

The aim of this study is to develop similar ideas for the the trigonometric functions which 
have singularities, namely, the   $ \tan, \cot, \sec $  and  $ \cosec $  functions.  
Our attention was drawn to this subject in an attempt to refine the 
inequality of Becker and Stark  \cite{Becker-Stark}
\begin{equation}						\label{eq:BeckerStark} 
\frac{8 / \pi^2} { 1 - x^2 }  \le    \frac {\tan (\pi x / 2) } { \pi x / 2 } 
		  \le   \frac{1} { 1 - x^2 } ,		\qquad  -1 < x < 1 .
\end{equation} 
Our main result (Theorem 1) is that the function  $  \dfrac {\tan (\pi x / 2) } { \pi x / 2 } $  
admits a sequence of upper and lower bounds and these inequalities are closely related to 
a Laurent series expansion in powers of  $ 1 - x^2 $.   The same approach will be applied 
also to  $ \cot, \sec $  and  $ \cosec $.

In spite of the similarity, there are several essential differences between our result and those  
of Ling Zhu  and  for  $ \sin x / x $  and of Baricz for the Bessel functions.	\\ 
(a)  Bessel functions and, consequently  $ \sin x = x^{-1/2} J_{1/2}(x) $,  are solutions 
of simple linear second order differential equations.  This is a useful tool in \cite{Baricz}.  
Unfortunately, we do not know any simple linear differential equation for the  tan  function.
On the contrary, the results for  $ \tan (\pi x / 2) / ( \pi x / 2 ) $  are influenced by 
its poles at $ x = \pm 1 $.  	\\
(b)  We start our study with the infinite series and the upper and lower bounds appear as consequences 
of this alternating series.  We do not use the monotone l'Hospital rule.  Our method uses partial 
fractions expansions.		\\

In analogy with  (\ref{eq:J-r})  and  (\ref{eq:Baricz})  we also expand 
  $ \dfrac {\tan (\pi x / 2) } { \pi x / 2 } $  into a Laurent series of  $ r^2 - x^2 $,
with  $ 0 < r < 1 $.  This expansion has different properties from that which corresponds
 to  $ r=1 $.


\section{ Inequalities and series for the  {\it tan}  function}
In this section we present an infinite series expansion for the tan function and
two sets of upper and lower bounds, one with an even number of terms and one with 
odd number of terms.

\begin{theorem}     
						\label{thm:tan}
(a) \ We have the infinite expansion 
\begin{equation}						\label{eq:tan-laurent} 
\frac {\tan (\pi x / 2) } { \pi x / 2 } 
	= \frac {8}{\pi^2} \left[  \frac{1} { 1 - x^2 } 
	  +  \sum_{ k=0 }^\infty  (-1)^k   \frac { T_{k+1} }{ 4^{k+1} } (1-x^2)^k \right] , 
\end{equation}
where 
\begin{equation}						\label{eq:Tp-zeta}
 T_p = (-1)^p \left[ 
	2\sum_{ m=1 }^{ [p/2] }  \binom{2p-2m-1}{p-1} \zeta(2m) - \binom{2p-1}{p-1} 
		\right],	\qquad 	p=1,2,\ldots \  . 
\end{equation}
All the coefficients  $ T_p $  are positive.	The infinite series in   (\ref{eq:tan-laurent}) 
converges for  $ |x| < 3 $.	 \rule[-3 mm]{0mm}{8mm}      \\
(b) \ For all even number  $ 2N $  and  $ -1 < x < 1 $,  we have the pair of inequalities
\begin{equation}						\label{eq:tan-2N} 
\begin{aligned} 
 \frac {8}{\pi^2}
 & \left[  \frac{1} { 1 - x^2 } 
	  +  \sum_{ k=0 }^{2N-1}  (-1)^k   \frac { T_{k+1} }{ 4^{k+1} } (1-x^2)^k 
	  + H_{2N+1} (1-x^2)^{2N}  \right]	\qquad\qquad\qquad\qquad	\\
	&  \qquad\qquad\qquad\qquad   
\le  \   \frac {\tan (\pi x / 2) } { \pi x / 2 }	\
\le  \   \frac {8}{\pi^2} \left[  \frac{1} { 1 - x^2 } 
	  +  \sum_{ k=0 }^{2N}  (-1)^k   \frac { T_{k+1} }{ 4^{k+1} } (1-x^2)^k 
	    \right] ,
\end{aligned}
\end{equation} 
where  
\begin{equation*}
	H_{2N+1}   = \frac{\pi^2}{8} 
		- \left(1 + \frac{T_1}{4} - \frac{T_2}{4^2} +- \cdots 
		- \frac{ T_{2N} }{ 4^{2N} } \right)  
\end{equation*}    
and   $ \displaystyle   0 < H_{2N+1} < \frac{ T_{2N+1} }{ 4^{2N+1} } $ .   
The left hand side inequality of  (\ref{eq:tan-2N})  is strict for  $ x=0 $.	\\
(c) \	 		\rule{0mm}{8mm}  
For all odd number  $ 2N+1 $  and  $ -1 < x < 1 $,  we have 
\begin{equation}				\label{eq:tan-2N+1}    
\begin{aligned} 
 \frac {8}{\pi^2} 
 &  \left[  \frac{1} { 1 - x^2 } 
	  +  \sum_{ k=0 }^{2N+1}  (-1)^k   \frac { T_{k+1} }{ 4^{k+1} } (1-x^2)^k 
	   \right]	\
\le   \   \frac {\tan (\pi x / 2) } { \pi x / 2 }	   \qquad\qquad		\\
	&  \qquad\qquad  
\le   \   \frac {8}{\pi^2} \left[  \frac{1} { 1 - x^2 } 
	  +  \sum_{ k=0 }^{2N}  (-1)^k   \frac { T_{k+1} }{ 4^{k+1} } (1-x^2)^k 
	  - H_{2N+2} (1-x^2)^{2N+1}   \right] ,
\end{aligned}
\end{equation} 
where  
$$
  H_{2N+2} = -\frac{\pi^2}{8} 
	+ \left(1 + \frac{T_1}{4} - \frac{T_2}{4^2} +- \cdots 
	+ \frac{ T_{2N+1} }{ 4^{2N+1} } \right) 
$$ 
and  $ \displaystyle  0 < H_{2N+2} < \frac{ T_{2N+2} }{ 4^{2N+2} } $ .   
The right hand side inequality of  (\ref{eq:tan-2N+1})  is strict for  $ x=0 $.
\end{theorem}

\noindent
{\bf  Examples.  \  }
The first coefficients are  $ T_1 = 1 $,  \ 
$ \displaystyle  T_2 = 2\zeta(2) - \binom{3}{1} = \pi^2 / 3 - 3 $  and   
$ \displaystyle  T_3 = -2 \binom{3}{2} \zeta(2) + \binom{5}{2} = -6\zeta(2) + 10 = 10 - \pi^2  $. 
Accordingly, 
\begin{align*}
&  H_1 = \frac{ \pi^2 } {8} - 1  \  ,			\qquad
   H_2 = -\left( \frac{ \pi^2 } {8} - 1 - \frac{T_1}{4}  \right)
	    = \frac{ 10 - \pi^2 } {8}	\   .
\end{align*}
The first two special cases of  (\ref{eq:tan-2N}) (for  $ 2N=0 $)  and  (\ref{eq:tan-2N+1}) 
(for  $ 2N+1 = 1 $)  are 
\begin{equation}					\label{eq:tan-ineq-0}
 \frac {8 / \pi^2} { 1 - x^2 } +  \frac{ \pi^2 - 8 } {\pi^2 } 
\le 	\frac {\tan (\pi x / 2) } { \pi x / 2 }	
\le  \frac {8/ \pi^2} { 1 - x^2 } + \frac { 2 }{  \pi^2 }   \ , \qquad
\end{equation}
  and
\begin{equation}					\label{eq:tan-ineq-1}
 \frac {8 / \pi^2} { 1 - x^2 } +  \frac{2}{ \pi^2 } 
	- \frac{ \pi^2 - 9 }{6 \pi^2} \, (1 - x^2)
\le 	\frac {\tan (\pi x / 2) } { \pi x / 2 }	
\le  \frac {8/ \pi^2} { 1 - x^2 } + \frac { 2 }{  \pi^2 }  
	-\frac{ 10 - \pi^2 } { \pi^2 } \, (1 - x^2)    . 
\end{equation}

\noindent
{\bf Proof of Theorem   \ref{thm:tan}  } 
For sake of convenience we divide the proof into several steps.

\noindent  
{\it (i) \  The infinite expansion for tan . }  \ 
Our starting point is the partial fraction expansion 
$  \displaystyle  \tan z = \sum_{ n=0 }^\infty \frac { 8z } { (2n+1)^2 \pi^2 - 4 z^2 } \  $ 
which holds in every domain which contains no poles.  It follows that 
\begin{equation}				\label{eq:tan-fractions} 
\frac {\tan (\pi x / 2) } { \pi x / 2 } 
	 =  \frac {8}{\pi^2}  \sum_{ n=0 }^\infty \frac {1} { (2n+1)^2 - x^2} 
	 =  \frac {8}{\pi^2} \left[  \frac{1} { 1 - x^2 } 
	    +  \sum_{ n=1 }^\infty   \frac { 1 }{ 4 n (n+1) + (1-x^2) } \right] . 
\end{equation} 
Each term of the series may be expanded into a geometric series, provided that 
$ | 1 - x^2 |   <  4n(n+1) $  for $ n=1,2,\ldots $,  i.e., for  $ |x|<3 $. This yields 
\begin{equation}				\label{eq:tan (pi x/2)}    
\begin{aligned} 
\frac {\tan (\pi x / 2) } { \pi x / 2 }
	=  &    \frac {8}{\pi^2}  \left[  \frac{1} { 1 - x^2 }  
	      +  \sum_{ n=1 }^\infty  \left( \sum_{ k=0 }^\infty 
	        \frac { (-1)^k (1-x^2)^k } { (4 n (n+1))^{k+1} } \right) \right]   \\
	=  &  \frac {8}{\pi^2} \left[  \frac{1} { 1 - x^2 } 
	      +  \sum_{ k=0 }^\infty  
		\left( \sum_{ n=1 }^\infty \frac{1} { n^{k+1} (n+1)^{k+1} } \right) 
	   	  \frac { (-1)^k }{ 4^{k+1} } (1-x^2)^k  \right] . 
\end{aligned}
\end{equation}


 The expansion (\ref{eq:tan-laurent})  will follow if we show that the sums   
$  \displaystyle
  T_{p}  = \sum_{n=1}^\infty \frac{1}{ n^p (n+1)^p } \ 
$ 
are given  explicitly by  (\ref{eq:Tp-zeta})  for  $ p=1,2,\ldots $.  The proof is based on 
the following identity   \cite[Ex. 39, p. 246]{Graham-Knuth-Patashnik}:  
{\sl  If  $ \ xy = x+y  $,  then for every}     $  \  m,k \ge 1 $, 
\begin{equation}						\label{eq:xm-yk}
    x^m y^k = \sum_{j=1}^m \binom{m+k-j-1}{k-1} x^j 
		+ \sum_{j=1}^k \binom{m+k-j-1}{m-1} y^j \ .
\end{equation}
(\ref{eq:xm-yk})  may be verified by mathematical induction on  $ m $  and  $ k $.  
We take in   (\ref{eq:xm-yk})    $ m=k=p $, 
   $  \displaystyle   x = -\frac{1}{n} $,  $ \displaystyle   y = \frac{1}{n+1} $,  
  and conclude that
\begin{equation}						\label{eq:x=-1/n}
  \frac{ (-1)^p }{ n^p (n+1)^p } 
	=  \sum_{j=1}^p   \binom{2p-j-1 }{ p-1 } 
		\left[	 \frac{ (-1)^j }{ n^j } + \frac{ 1 }{ (n+1)^j }    \right] .
\end{equation}
To calculate  $ T_p $  we sum  (\ref{eq:x=-1/n})  for  $ n=1, 2, \ldots $ and  
separate  $ j=1 $  from  $ j= 2,\ldots, p $ :
\begin{align*}	
  (-1)^p T_{p}  
&  = (-1)^p 
\sum_{n=1}^\infty \sum_{j=1}^p   \binom{2p-j-1 }{ p-1 } 
		\left[ \frac{ (-1)^j }{ n^j } + \frac{ 1 }{ (n+1)^j } \right]	\\
&  =   \binom{2p-2 }{ p-1 } \sum_{n=1}^\infty 
				\left[  -\frac{1}{n} + \frac{1}{n+1}  \right]
	 +  \sum_{j=2}^p   \binom{2p-j-1 }{ p-1 } 
			\Big[ 	(-1)^{j} \zeta(j) + \big( \zeta(j) - 1 \big) \Big]   \\
&  =  \binom{2p-2}{p-1} \cdot (-1) 
	+ \sum_{ j = 2 }^{ p }  \binom{2p-j-1}{p-1} 
	   \Big[ \big( (-1)^{j} +1 \big)  \zeta(j) - 1 \big) \Big]   	\\
& = 2 \sum_{ 2 \le 2m \le p } \binom{2p-2m-1}{p-1} \zeta(2m) 
	- \sum_{j=1}^{p}  \binom{2p-j-1}{p-1}  . 
\end{align*}
Since 
$  \  \displaystyle \sum_{j=1}^{p}  \binom{2p-j-1}{p-1} = \binom{2p-1}{p-1} $,  \  
 formula  (\ref{eq:Tp-zeta})  is verified and  part (a) of the theorem follows.

We remark that similar sums   
$ \displaystyle  \sum_{n=1}^\infty \frac{1}{ (n - 1/2 )^p (n + 1/2 )^p } $ 
are mentioned in  \cite[p. 70, (373)]{Jolley}   for  $ p=1,2,3 $ and $ 4 $.

\noindent  
{\it (ii) \ The inequalities . }  \ 
Let us write the series (\ref{eq:tan-laurent}) as 
\begin{equation*}					
\frac {\tan (\pi x / 2) } { \pi x / 2 } 
 = \frac {8}{\pi^2} \left[  \frac{1} { 1 - x^2 } 
	  +  \sum_{ k=0 }^m  (-1)^k   \frac { T_{k+1} }{ 4^{k+1} } (1-x^2)^k 
	  + R_{m+1}(x)  \right] , 
\end{equation*}
with the remainder term 
$$
  R_{m+1}(x) = (-1)^{m+1} (1-x^2)^{m+1} 
    \sum_{ \ell=0 }^\infty  (-1)^\ell \frac { T_{\ell+m+2} }{ 4^{\ell+m+2} } (1-x^2)^\ell \ . 
$$
For each  $ x $,  $ 0 \le x \le 1 $,  the series 
$$
  F_{m+1}(x) = \sum_{ \ell=0 }^\infty  (-1)^\ell \frac{ T_{\ell+m+2} }{ 4^{\ell+m+2} } (1-x^2)^\ell 
$$
and its derivative series 
$$
  F'_{m+1}(x) = 2x \sum_{ \ell=1 }^\infty  (-1)^{\ell-1} \ell \frac{ T_{\ell+m+2} }{ 4^{\ell+m+2} } 
		(1-x^2)^{\ell-1} 
$$
are alternating series with smaller and smaller terms.   Indeed, 
$$
  \frac{1}{2^p} < T_p = \sum_{ n=1 }^\infty \frac{1} { n^p (n+1)^p } 
  		< \frac{1}{2^p} \zeta(p) < \frac{1}{2^p} \zeta(2) < \frac{1}{ 2^{p-1} } ,
$$
thus 
$$
 \frac{ T_{\ell+m+2} }{ 4^{\ell+m+2} } > \frac{1}{ 2^{\ell+m+2} \, 4^{\ell+m+2} } 
      >  \frac{1}{ 2^{\ell+m+2} \, 4^{\ell+m+3} } > \frac{ T_{\ell+m+3} }{ 4^{\ell+m+3} } 
$$
and similarly
$$
 \ell \frac{ T_{\ell+m+2} }{ 4^{\ell+m+2} } > (\ell+1) \frac{ T_{\ell+m+3} }{ 4^{\ell+m+3} } .
$$
Thus, by properties of alternating series,  $ F_{m+1}(x) $    and   $ F'_{m+1}(x) $ 
are positive functions for  $ 0 \le x \le 1 $.   Consequently    $ F_{m+1}(x) $  increases 
for  $ 0 \le x \le 1 $,   and is bounded from below and from above by  
$$
	0 <  \sum_{ \ell=0 }^\infty  (-1)^\ell \frac { T_{\ell+m+2} }{ 4^{\ell+m+2} }  \ 
	  =  F_{m+1}(0) \le  F_{m+1}(x)   \le  F_{m+1}(1)  =     \frac{ T_{m+2} } { 4^{m+2} }  \ .
$$
The lower bound of  $ F_{m+1}(x) $,   \
$  \frac{ T_{m+2} } { 4^{m+2} } - \frac{ T_{m+3} } { 4^{m+3} } +- \cdots $ 
may be written as a finite sum:  For  $ x=0 $, (\ref{eq:tan-laurent}) 
yields 
\begin{equation}					\label{eq:Tk-sum}
  1 = \frac {8}{\pi^2}\left[ 1 + \sum_{ k=0 }^\infty  (-1)^k  \frac { T_{k+1} }{ 4^{k+1} } \right] ,
\end{equation}
so
\begin{equation}					\label{eq:H(m+2)}
\begin{aligned}	
  0 <  H_{m+2}  
&   := \frac{ T_{m+2} } { 4^{m+2} } - \frac{ T_{m+3} } { 4^{m+3} } +- \cdots \\
&   = (-1)^{m+1}   \left[ \frac {\pi^2}{8} 
		-\left( 1 + \frac{ T_1 } {4} - \frac{ T_2 } { 4^2 } +- \cdots	
		    + (-1)^m \frac{ T_{m+1} } { 4^{m+1} } 
		\right)  \right] . 
\end{aligned}
\end{equation}
Consequently  $  \  R_{m+1}(x) = (-1)^{m+1} (1-x^2)^{m+1} F_{m+1}(x)  \  $  
is bounded for  $ m=2N $  and for  $ m = 2N-1 $,  respectively, by
\begin{equation*}
\begin{aligned}
  -(1-x^2)^{2N+1} \frac{ T_{2N+2} } { 4^{2N+2} }  &	\  \le  \  R_{2N+1}(x)  
		\  \le  \   -(1-x^2)^{2N+1} H_{2N+2} < 0 ,	\\
  0 < (1-x^2)^{2N}  H_{2N+1}  &		\  \,  \le  \  \  R_{2N}(x)  \ \
  \le   \  (1-x^2)^{2N}  \frac{ T_{2N+1} } { 4^{2N+1} }    \ .
\end{aligned}
\end{equation*} 
and parts (b) and (c) of the theorem follow.


The remainder term  satisfies 
$
	| R_{m+1}(x) | \le (1-x^2)^{m+1} \dfrac{ T_{m+2} } { 4^{m+2} } 
			\le   \dfrac{ (1-x^2)^{m+1} } { 2^{3m+5} } ,
$  
hence it decreases rapidly, in particular near the poles  $ x= \pm 1 $.  For example,  
for $ x=0.9 $  and  $ m=4 $  we have  
$ R_5(0.9) < 0.19^5 / 2^{17} \approx 2 \cdot 10^{-9} $.

The bounds in  (\ref{eq:tan-2N+1})  and  (\ref{eq:tan-2N})  are also very precise
near the poles  $ x = \pm 1 $.  Indeed, the difference between the  $ m $-th  upper 
and lower bounds  ($m$  odd or even) is less then 
$ \  
 \dfrac{ T_{m+3} } { 4^{m+3} } (1-x^2)^{ m+1 } < \dfrac{ 1 } { 2^{3m+8} } (1-x^2)^{ m+1 } $ .

It is interestiong to compare the convergence of the series  (\ref{eq:tan-laurent})  with
that of the Taylor series   \cite[\ Section \ 3:14]{Sp-Ol} 
\begin{equation}						\label{eq:tan-lambda}
\frac {\tan (\pi x / 2) } { \pi x / 2 } 
    = 1 + \frac {8}{\pi^2} \sum_{ k=1 }^\infty \lambda(2k+2) x^{2k} , 
	\qquad    |x| < 1 , 
\end{equation}
where  $ \lambda(p) = \sum_{ n=0 }^\infty  1/(2n+1)^{p} $  denotes the Dirichlet lambda 
function.   Since  $ 1 < \lambda(p) < 1 + 2/3^p $,  the remainder of the series  
(\ref{eq:tan-lambda})  after  $m$  terms is 
$$
\frac {8}{\pi^2} \sum_{ k=m+1 }^\infty \lambda(2k+2) x^{2k} 
= \frac {8}{\pi^2} \frac{ x^{2m+2} }{ 1 - x^2 } \, ( 1 + o(1) ) .
$$
Here  $ o(1) $  is uniformly small for large values of  $m$.  On the other hand, we saw 
above that the remainder of the series  (\ref{eq:tan-laurent})  after  $m$  terms is 
$ \dfrac {8}{\pi^2} (1-x^2)^{ m+1 } \dfrac{ O(1) } { 2^{3m} } $.

After extracting  $m$-th root from the two remainders for large  $m$, we compare, up to  $ 1+o(1) $,  
$ x^2 $  with  $ (1 - x^2)/8 $.  The consequence is that after  $m $  terms, for large  
$m$,  the Taylor series  (\ref{eq:tan-lambda})  yields a smaller remainder for  $ |x| < 1/3 $,
while for  $ 1/3 < |x| < 1 $,  the series  (\ref{eq:tan-laurent})  has a smaller remainder.
\qed



Due to  (\ref{eq:J-r})  and   (\ref{eq:Baricz})  it is also natural to expand  
$  \  \dfrac {\tan (\pi x / 2) } { \pi x / 2 }  \ $  into  power series of  $ r^2 - x^2 $.
Formally we follow  (\ref{eq:tan-fractions})  and  (\ref{eq:tan (pi x/2)}):

\begin{equation}				\label{eq:tan-r2-x2)}  
\begin{aligned} 
\frac {\tan (\pi x / 2) } { \pi x / 2 } 
	 =  & \frac {8}{\pi^2}  \sum_{ n=0 }^\infty \frac {1} { (2n+1)^2 - x^2} 
= \frac {8}{\pi^2}  \sum_{ n=0 }^\infty \frac {1} { \big((2n+1)^2 - r^2 \big) + ( r^2 - x^2 ) } \\
=  &  \frac {8}{\pi^2} 
		\left[   \sum_{ k=0 }^\infty  
			\left( \sum_{ n=0 }^\infty 
			\frac{1} { ( n + \frac{1-r}{2} )^{k+1} (n+ \frac{1+r}{2} )^{k+1} } 
			\right) 
	   	  	\frac { (-1)^k }{ 4^{k+1} } (r^2-x^2)^k  
		\right] . 		\\
=  &   \frac {8}{\pi^2}  
	    \sum_{ k=0 }^\infty  (-1)^k   \frac { \tilde{T}_{k+1}(r) }{ 4^{k+1} } (r^2-x^2)^k   \,  ,
\end{aligned}
\end{equation} 
with
\begin{equation}				\label{eq:T-p-r}    
 	\tilde{T}_{p}(r) 
	= \sum_{ n=0 }^\infty \frac{1} { ( n + \frac{1-r}{2} )^{p} (n+ \frac{1+r}{2} )^{p} } .
\end{equation} 
Note that  $ \tilde{T}_{p}(r) $  differs from  $ T_p $  by it's range of summation and it is
not defined for  $ r=1 $.

Equations  (\ref{eq:tan-r2-x2)})  and  (\ref{eq:tan (pi x/2)})  are superficially similar.
However, they describe different phenomena.  (\ref{eq:tan (pi x/2)})  describes the function 
$ \tan (\pi x / 2) / (\pi x / 2) $  in the interval  $ (-3, 3) $, including its singularities at 
$ x = \pm 1 $.  Expansion  (\ref{eq:tan-r2-x2)})  is valid when  $ |r^2 - x^2 | < (2n+1)^2 - r^2 $
for all  $ n=0,1,2, \ldots $ . Assuming  $ r < 1 $,  this holds when  $ |r^2 - x^2 | < 1 - r^2 $,
i.e., 
		$$ 	2 r^2 - 1 < x^2 < 1 ,	 $$ 
which  includes the points $ x = \pm r $  and their neighborhoods.  Note that in the complex plane 
$ |r^2 - z^2 | < C $  is the interior of a lemniscate with focal points  $ z = \pm r $.

Inequalities analogous to  (\ref{eq:tan-2N})  and  (\ref{eq:tan-2N+1})  may be deduced from
expansion  (\ref{eq:tan-r2-x2)}).  We do not attempt to calculate the sums  (\ref{eq:T-p-r})  
explicitly as functions of  $ r $.  Rather we shall obtain a recursive formula for the coefficients 
$ \tilde{T}_{p}(r) $.   

The function  $ \tan x $ satisfies the first order, quadratic Riccati differential equation
$ y' = 1 + y^2 $.  It easy to see that  $ u(x) = \dfrac {\tan (\pi x / 2) } { \pi x / 2 } $ 
satisfies  $ x \dfrac{du}{dx} = 1 - u + \dfrac{\pi^2}{4} x^2 u^2 $  for  $ |x|<1 $.  
If we let  $ z = r^2 - x^2 $  and  $ u(x) = v(z) $,  then  
$ v(z)  $
satisfies 
$$
	-2 (r^2 - z) \frac{dv}{dz} = 1 - v + \dfrac{\pi^2}{4} (r^2 - z) v^2   
$$ 
for  $ -1 + r^2 < z < r^2 $.  Substitution of    
$ v(z) = \dfrac {8}{\pi^2} \sum_{ k=0 }^\infty  (-1)^k  \dfrac { \tilde{T}_{k+1}(r) }{ 4^{k+1} } \, z^k$ ,  
$ v(0) = \displaystyle  \frac {\tan (\pi r / 2) } { \pi r / 2 } $
and comparison of coefficients lead to the following quadratic recursion formula for  
$ \tilde{T}_k(r) $,  \   $ 0 < r < 1 $:
\begin{equation}				\label{eq:Tr)}  
\begin{aligned} 
	     \tilde{T}_1(r) & = \frac{\pi^2}{2} \, \frac {\tan (\pi r / 2) } { \pi r / 2 },	\\
	 r^2 \tilde{T}_2(r) & = \pi^2 - 2 \tilde{T}_1(r) + r^2 \tilde{T}_1^2(r),		\\
(k+1)r^2 \tilde{T}_{k+2}(r) & = - (4k+2) \tilde{T}_{k+1} + r^2 \sum_{j=0}^k \tilde{T}_{j+1} \tilde{T}_{k-j+1} 
					+ 4 \sum_{j=0}^{k-1} \tilde{T}_{j+1} \tilde{T}_{k-j} ,
		\quad   k=1,2,\ldots  	\hspace{-10 mm} 
\end{aligned}
\end{equation} 
%


\section{Inequalities for the secant function. }

Expansions and inequalities of similar type are available for the  $ \cot, \sec $  and  
$ \cosec $  functions.  In this section we present the results for the  secant  function.
We remind, in analogy with inequalities  (\ref{eq:BeckerStark})  and  (\ref{eq:Chen-Qi-cot} ),  
that
\begin{equation}					\label{eq:sec-BeckerStark} 
\frac{1} { 1 - x^2 }  \le    \sec \frac{\pi x }{ 2 }   \le   \frac{4 / \pi} { 1 - x^2 } ,		
		\qquad  -1 < x < 1 .
\end{equation}
To prove  (\ref{eq:sec-BeckerStark}), it is sufficient to show that the even function 
 $ ( 1 - x^2 ) \sec (\pi x / 2) $   increases for  $ 0 \le x \le 1 $.  This indeed holds,  
since by  (\ref{eq:BeckerStark}) 
\begin{equation}					\label{eq:sec-derivative} 
  \frac{d}{dx} \left( (1-x^2)  \sec(\pi x / 2)   \right) 
= 	2 x \sec(\pi x / 2) 
	\left[ \frac{ \pi^2}{8 } (1-x^2) \frac {\tan (\pi x / 2) } { \pi x / 2 } - 1 \right] 
		\ge 0 .
\end{equation} 
The constants in   (\ref{eq:sec-BeckerStark})  follow now by the limit values for 
$ x=0 $  and as  $ x \to 1^- $.

Our aim is to generalize the right hand side of  (\ref{eq:sec-BeckerStark}) 
into an infinite sequence of inequalities.

\begin{theorem}   
								\label{thm:sec}
(a) \ We have the infinite expansion 
\begin{equation}						\label{eq:sec-laurent} 
 \sec (\pi x / 2) 
	= \frac{4}{\pi} \left[  \frac{1} { 1 - x^2 } 
	  -  \sum_{ k=0 }^\infty  (-1)^k   \frac { S_{k+1} }{ 4^{k+1} } (1-x^2)^k \right] , 
\end{equation}
where  
\begin{equation}						\label{eq:Sp-eta} 
 S_p 
= (-1)^p \left[ \binom{2p-2 }{ p-2 } - \binom{2p-2 }{ p-1 } 
		  + 2 \sum_{ 2 \le 2m \le p }  
		   \Big[  \binom{2p-2m-2 }{ p-2 } - \binom{2p-2m-2 }{ p-1 } \Big] \eta(2m) 
	 \right]  
\end{equation} 
for  $ p = 1,2, \ldots $  and  $ \eta $ denotes the Dirichlet eta function.
All the coefficients  $ S_p $  are positive.		\rule[-3 mm]{0mm}{8mm}	  \\
%
(b) \ For all even number  $ 2N $  and  $ \  -1 < x < 1 $,  we have the pair of inequalities
\begin{equation}						\label{eq:sec-2N} 
\begin{aligned} 
 \frac{4}{\pi} 
 & \left[ \frac{1} { 1 - x^2 } 
	  -  \sum_{ k=0 }^{2N}  (-1)^k  \frac { S_{k+1} }{ 4^{k+1} } (1-x^2)^k
   \right]	   
\le  \    \sec (\pi x /2 ) 	\\
& \qquad\qquad   
 \le  \   \frac{4}{\pi} \left[  \frac{1} { 1 - x^2 } 
	  -  \sum_{ k=0 }^{2N-1}  (-1)^k   \frac { S_{k+1} }{ 4^{k+1} } (1-x^2)^k 
	  - J_{2N+1} (1-x^2)^{2N}
	    \right] ,
\end{aligned}
\end{equation} 
where  
\begin{equation*}
	J_{2N+1}   = -\frac{\pi}{4} 
				+ \left( 1 - \frac{S_1}{4} + \frac{S_2}{4^2} -+ \cdots 
					+ \frac{ S_{2N} }{ 4^{2N} } 
 		  \right)  
\end{equation*}    
and   $ \displaystyle   0 < J_{2N+1} < \frac{ S_{2N+1} }{ 4^{2N+1} } $ .   
The right hand side inequality of  (\ref{eq:sec-2N})  is strict for  $ x=0 $.	\\
(c) \	 		\rule{0mm}{8mm}  
For all odd number  $ 2N+1 $  and  $  \  -1 < x < 1 $,  we have 
\begin{equation}				\label{eq:sec-2N+1}    
\begin{aligned} 
 \frac{4}{\pi} 
 &  \left[ \frac{1} { 1 - x^2 } 
	  -  \sum_{ k=0 }^{2N}  (-1)^k   \frac { S_{k+1} }{ 4^{k+1} } (1-x^2)^k 
	  + J_{2N+2} (1-x^2)^{2N+1}  
   \right]							\\
&  \qquad\qquad   
\le   \    \sec (\pi x /2 ) 	   
\le   \   \frac{4}{\pi} \left[  \frac{1} { 1 - x^2 } 
	  -  \sum_{ k=0 }^{2N+1}  (-1)^k \frac { S_{k+1} }{ 4^{k+1} } (1-x^2)^k 
    \right] ,
\end{aligned}
\end{equation} 
where  
$$
  J_{2N+2} = \frac{\pi}{4} 
		- \left(1 - \frac{S_1}{4} + \frac{S_2}{4^2} -+ \cdots 
			- \frac{ S_{2N+1} }{ 4^{2N+1} } 
	  \right) 
$$ 
and  $ \displaystyle  0 < J_{2N+2} < \frac{ S_{2N+2} }{ 4^{2N+2} } $ .   
The left hand side inequality of  (\ref{eq:sec-2N+1})  is strict for  $ x=0 $.
\end{theorem}

\noindent
{\bf  Examples.  \  }
The first coefficients are  $ S_1 = 1 $,  \ 
$ \displaystyle  S_2 = -1 + 2\eta(2) = \pi^2 / 6 - 1 $  and   
$ \displaystyle  S_3 = 2 -2 \eta(2) = 2 - \pi^2 / 6 $. 
Accordingly, 
\begin{align*}
&  J_1 = -\frac{ \pi }{4} + 1  \  ,			\qquad
   J_2 = \frac{ \pi }{4} - 1 + \frac{1}{4} = \frac{ \pi - 3 }{4}	\   .
\end{align*}
For  $ 2N=0 $  we get on  $  \  -1 < x < 1 $ 
$$
  \frac {4}{\pi} \left[  \frac{1} { 1 - x^2 } - \frac{1}{4}  \right] 
	\le	\sec (\pi x / 2) 
	\le	\frac {4}{\pi} \left[  \frac{1} { 1 - x^2 } - 1 + \frac{\pi}{4} \right]  
$$ 
and for  $ 2N+1=1$, 
$$
  \frac {4}{\pi} \left[  \frac{1} { 1 - x^2 } - \frac{1}{4} + \frac{ \pi - 3 }{4}(1 - x^2) \right] 
	\le	\sec (\pi x / 2) 
	\le	\frac {4}{\pi} \left[  \frac{1} { 1 - x^2 } - \frac{1}{4} 
		+ \frac{\pi^2 -6}{96} (1 - x^2)  \right]  \  .
$$

\noindent
{\bf Proof.}    We start with the partial fraction expansion 
$  \displaystyle
	\sec z = \pi \sum_{ n=0 }^\infty \frac { (-1)^n (2n+1) } { (n+1/2)^2 \pi^2 -  z^2 } $. 
Then for  $  \  -1 < x < 1 $, 
\begin{equation*}			
\begin{aligned} 
\sec (\pi x / 2)  
	& =  \frac {4}{\pi}  \sum_{ n=0 }^\infty \frac {(-1)^n (2n+1)} { (2n+1)^2 - x^2}	\\
	& = \frac {4}{\pi} \left[  \frac{1} { 1 - x^2 } 
	    - \sum_{ n=1 }^\infty   \frac { (-1)^{n-1} (2n+1) }{ 4 n (n+1) + (1-x^2) } \right] 	\\ 
	 & = \frac {4}{\pi} \left[  \frac{1} { 1 - x^2 } 
	    - \sum_{ k=0 }^\infty  
		\left( \sum_{ n=1 }^\infty \frac{ (-1)^{n-1} (2n+1) } { n^{k+1} (n+1)^{k+1} } \right) 
	   	(-1)^k   \frac { (1-x^2)^k }{ 4^{k+1} } \right] 				\\
 	& = \frac {4}{\pi} \left[  \frac{1} { 1 - x^2 } 
	    -  \sum_{ k=0 }^\infty  (-1)^k   \frac { S_{k+1} }{ 4^{k+1} } (1-x^2)^k \right] , 
\end{aligned}
\end{equation*} 
  where  
$  \  \displaystyle  
     S_p = \sum_{ n=1 }^\infty \frac{ (-1)^{n-1} (2n+1) } { n^p (n+1)^p } > 0  \  $ 
for $ p=1,2, \ldots $ .

The coefficients  $ S_p $  may be presented explicitly by finite sums of the Dirichlet 
$ \eta $  function.  Let us multiply   (\ref{eq:x=-1/n})  by  $ 2n+1 \equiv 2(n+1) - 1 $ 
and separate  $ j=1 $ from  $ 2 \le j \le p $:
\begin{equation*}						
\begin{aligned} 
  \frac{ (-1)^p (2n+1)}{ n^p (n+1)^p } 
	= &  - \binom{2p-2 }{ p-1 } \left[ \frac{ 1 }{ n } + \frac{ 1 }{ n+1 } \right] 	\\
	  & +  \sum_{j=2}^p   \binom{2p-j-1 }{ p-1 } 
		\left[	2\frac{ (-1)^j }{ n^{j-1} } + \frac{ (-1)^j }{ n^j } 
			+ 2\frac{ 1 }{ (n+1)^{j-1} } - \frac{ 1 }{ (n+1)^j }   \right] .
\end{aligned}
\end{equation*}
We multiply the last equation by  $ (-1)^{n-1} $  and sum it for  $ n = 1,2, \ldots $ .
With the Dirichlet eta function  
\begin{equation}							\label{eq:eta}
  	\eta (p) =  \sum_{ n=1 }^\infty  \frac { (-1)^{n-1} } { n^p },  \qquad
	  \sum_{ n=1 }^\infty  \frac { (-1)^{n-1} } { (n+1)^p } = -\eta (p) +1 ,
	\qquad  p \ge 1  , 
\end{equation} 
 we get
\begin{equation*}
\begin{aligned} 
(-1)^p & S_p  = - \binom{2p-2 }{ p-1 } \left[ \eta(1) + \big( -\eta(1) + 1 \big)  \right] 	\\
&	+  \sum_{ j=2 }^p  \binom{2p-j-1 }{ p-1 } 
		\Big[   2(-1)^j \eta(j-1) + (-1)^j \eta(j) 
		      + 2 \big(-\eta(j-1) + 1 \big) - \big(-\eta(j) + 1 \big)   \Big]	\\ 
&  =  -\binom{2p-2 }{ p-1 }  +  \sum_{ j=2 }^p  \binom{2p-j-1 }{ p-1 } 
	+ \sum_{ j=2 }^p  \binom{2p-j-1 }{ p-1 } 
	     \Big[  2 \big( (-1)^j - 1 \big)  \eta(j-1) + \big( (-1)^j + 1 \big) \eta(j)     \Big]	\\ 
&  =  - \binom{2p-2 }{ p-1 } + \binom{2p-2 }{ p }  
	  +  \sum_{ 2 \le 2m \le p }  
	       2 \Big[  \binom{2p-2m-1 }{ p-1 } - 2 \binom{2p-2m-2 }{ p-1 } \Big] \eta(2m)  \\
&  =  \binom{2p-2 }{ p-2 } - \binom{2p-2 }{ p-1 }   
	+ 2 \sum_{ 2 \le 2m \le p }  
	        \Big[  \binom{2p-2m-2 }{ p-2 } - \binom{2p-2m-2 }{ p-1 } \Big] \eta(2m) .	       
\end{aligned}
\end{equation*} 
Here we used the facts that    
$   \displaystyle \sum_{ j=2 }^p  \binom{2p-j-1 }{ p-1 } = \binom{2p-2 }{ p } = \binom{2p-2 }{ p-2 }$ 
and that 
$$   	\binom{2p-2m-1 }{ p-1 } - 2 \binom{2p-2m-2 }{ p-1 } 
      = \binom{2p-2m-2 }{ p-2 } - \binom{2p-2m-2 }{ p-1 } .
$$
Thus,  (\ref{eq:Sp-eta})  is proved.

The inequalities  (\ref{eq:sec-2N})  and   (\ref{eq:sec-2N+1})  are verified as in 
Theorem  \ref{thm:tan}. 
\qed

The coefficients  $ T_p $  and  $ S_p $  are easily related.  In the (\ref{eq:sec-derivative}) 
we substitute  $ \sec( \pi x/2) $  and  $ \dfrac{ \tan( \pi x / 2 ) } { \pi x / 2 } $, 
respectively, by the series  (\ref{eq:tan-laurent})  and   (\ref{eq:sec-laurent})  and 
compare powers of  $ 1 - x^2 $,  the result is 
$$
   S_{n+1} = T_{n+1} + \sum_{k=0}^{n-1} T_{k+1} S_{n-k}  \  .  
$$


\section{Inequalities for other trigonometric functions. }

We outline some other results which may be verified by similar arguments. 
For the cotangent function, let us recall the inequality  \cite[Eq. (17)]{Chen-Qi-2004}
\begin{equation}				\label{eq:Chen-Qi-cot}    
	\frac {2 x^2 }{ 1-x^2 }  \le  1 - \pi x \cot (\pi x ) 	
				 \le \frac{\pi^2}{3} \frac {x^2 }{ 1-x^2 } \  .
\end{equation}
It is possible generalize the left hand side of  (\ref{eq:Chen-Qi-cot})  into the infinite 
expansion
%
\begin{equation}				\label{eq:cot-laurent}    
1 -  \pi x \cot (\pi x ) 
	=  {2} {x^2}  \left[  \frac{1} { 1 - x^2 } 
	  +  \sum_{ k=0 }^{\infty}  (-1)^k   \frac { C_{k+1} }{ 4^{k+1} } (1-x^2)^k 
	    \right]   \ ,
\end{equation}
\begin{equation}					\label{eq:Cp-zeta}
\begin{aligned} 
 C_p 
&   	= \sum_{ n=2 }^\infty \frac{ 2^{2p} } { (n^2 - 1)^p  } 
		= \sum_{ n=1 }^\infty \frac{ 2^{2p} } { n^p (n+2)^p  } 		\\
&	= (-1)^p \left[ 
		\sum_{ 2 \le 2m \le p } 2^{ 2m+1 } \binom{2p-2m-1}{p-1} \zeta(2m) 
		- \binom{2p-1}{p-1} - 2^{  2p-1 } 
		\right]	, \quad   p \ge 1 .
\end{aligned}
\end{equation}
All  $ C_p $-s  are positive. The infinite series in  (\ref{eq:cot-laurent})  converges 
for  $ |x| < 2 $.

For example,  $ C_1=3, \ C_2 = -11 + 8\zeta(2) = \frac{4}{3}\pi^2 - 11 , \
C_3 = 42 - 24\zeta(2) = 42 - 4\pi^2 $.

Inequalities in the style of Theorems \ref{thm:tan}  and  \ref{thm:sec}  may be proved for 
$ 1 -  \pi x \cot (\pi x ) $.


(\ref{eq:cot-laurent})  is obtained from the expansion 
$ \displaystyle 
    \pi x \cot (\pi x ) = 1 - 2x^2 \sum_{ n=1 }^\infty   \frac { 1 }{n^2-x^2 } 
$. 
For the calculation of  $ C_p $  in  (\ref{eq:Cp-zeta}),  it is useful to take in  (\ref{eq:xm-yk}),
  $ m=k=p $, 	$ \displaystyle   x = \frac{-2}{n} , y = \frac{2}{n+2} $   
and sum the result for  $ n=1,2,  \ldots $.
For the summation of the  series, one has to use the binomial identity  
%
$  \   \displaystyle  \sum_{j=1}^{p}  \binom{ 2p-j-1 }{p-1} 2^j = 2^{2p-1} , \  $  
\cite[ p. 167, Eq. 5.20]{Graham-Knuth-Patashnik}.
\qed


  For the  cosecant  function we remind, in analogy with 
  (\ref{eq:BeckerStark}),  (\ref{eq:sec-BeckerStark})  and  (\ref{eq:Chen-Qi-cot}), 
the double inequality  \cite[Eq. (19)]{Chen-Qi-2004} 
\begin{equation}					\label{eq:Chen-Qi-cosec} 
 \frac{ (\pi^2/6) x^2} { 1 - x^2 }  \le  \pi x  \cosec( \pi x ) - 1  
					\le   \frac{2 x^2} { 1 - x^2 } ,		
		\qquad  -1 < x < 1 .
\end{equation}
The right hand side of  (\ref{eq:Chen-Qi-cosec})  may be generalized into the expansion 
%
\begin{equation}						\label{eq:cosec-laurent} 
\pi x \cosec (\pi x ) - 1
	= 2 x^2 \left[  \frac{1} { 1 - x^2 } 
	      -  \sum_{ k=0 }^\infty  (-1)^k   \frac { D_{k+1} }{ 4^{k+1} } (1-x^2)^k \right] , 
\end{equation}
where  
\begin{equation}						\label{eq:Dp-eta}
\begin{aligned} 
 D_p 
&   =   \sum_{ n=2 }^\infty \frac{ (-1)^{n} 2^{2p} } { (n^2 - 1)^p  } 
		= \sum_{ n=1 }^\infty \frac{ (-1)^{n-1} 2^{2p} } { n^p (n+2)^p  } 	\\
&	= (-1)^p \left[ 
		\sum_{ 2 \le  2m \le p } 2^{2m+1 } \binom{2p-2m-1 }{ p-1 } \eta(2m)  
			- 2^{2p-1} + \binom{2p-1 }{ p-1 }
		\right] .  
\end{aligned}
\end{equation}

Here  $ D_1=1, \ D_2 = -5 + 8\eta(2) = \frac{2}{3}\pi^2 - 5 , \
D_3 = 22 - 24\eta(2) = 22 - 2\pi^2 $.


Expansion (\ref{eq:cosec-laurent})  follows from the partial fraction expansion
$  \displaystyle
	\cosec z = \frac{1}{z} - 2z \sum_{ n=1 }^\infty \frac { (-1)^n } { n^2 -  z^2 } $. 
$ D_p $  is calculated similarly to  $ C_p $. 
\qed

If we substitute the series  (\ref{eq:cot-laurent})  and  (\ref{eq:cosec-laurent}) 
into the identity 
$$
 x \frac{d}{dx} \left( (1 - x^2) \pi x \cosec(\pi x) \right) 
= \pi x \cosec(\pi x) \Big[ -2 x^2 + (1 - x^2)( 1 - \pi x \cot(\pi x) ) \Big]  
$$
 and compare powers of
$ 1-x^2 $,  it follows that  $ C_n $  and  $ D_n $  are related by 
$$
n D_n + 4n D_{n-1} 
   = C_n + 2 C_{n-1} + \sum_{k=1}^{n-1} D_k C_{n-k} +  \sum_{k=1}^{n-2} D_k C_{n-k-1} \ . 
$$


\section{Inequalities for Bessel functions}

Finally we show that our methods are applicable also to Bessel functions and they
provide a short proof of  (\ref{eq:Baricz}),  based on other principles than those
of Baricz and Wu in  \cite{Baricz-Wu}.

  An expansion of the form 
\begin{equation}				\label{r2-x2-expansion}
	x^{-p} J_p(x) = \sum_{k=0}^\infty  c_k (r^2 - x^2)^k 
\end{equation} 
holds for all  $ x $  in the complex plane.  To see this, put  $ z = r^2 - x^2 $,  
i.e.,  $ x = \sqrt{ r^2 - z } $.  Since  $ x^{-p} J_p(x) $ 
is an even analytic function in the whole complex plane, it follows that
	$ x^{-p} J_p(x) \Big|_{ x = \sqrt{ r^2 - z } } $  
is analytic for all  $ z $  and may be expanded as  $ \sum_{k=0}^\infty  c_k z^k $. 
Consequently  (\ref{r2-x2-expansion})  holds for all  $ x $.
To calculate the  $c_k$-s,  let us apply the identity  \cite[Eq. 9.1.30]{Abramowitz-Stegun}
$$
    \left( \frac{1}{x} \frac{d}{dx} \right)^m \left( x^{-p} J_p(x) \right)
=	(-1)^m x^{ -(p+m) } J_{ p+m }(x) , \qquad m = 1,2,\ldots 
$$
to expansion  (\ref{r2-x2-expansion}).  Since   
$ \left( \frac{1}{x} \frac{d}{dx} \right)^m  (r^2 - x^2)^k 
		= (-2)^m k(k-1) \cdots (k-m+1) (r^2 - x^2)^{k-m} $, 
we get
$$
  (-1)^m x^{ -(p+m) } J_{ p+m }(x) 
  = (-2)^m \sum_{ k=m }^\infty  c_k \, k(k-1) \cdots (k-m+1) (r^2 - x^2)^{k-m} . 
$$
For  $ x=r $  we obtain that  $ c_m = \dfrac{ r^{-(p+m) } J_{p+m}(r) } { 2^m m! } $. 
Consequently, 
\begin{equation}				\label{Baritz-expansion}
	x^{-p} J_p(x) = \sum_{k=0}^\infty \frac{ r^{-(p+k) } J_{p+k}(r) } { 2^k k! } (r^2 - x^2)^k ,
	\qquad  |x| < \infty .
\end{equation}

To prove the inequalities (\ref{eq:Baricz}), we have to estimate the remainder of the 
expansion  (\ref{Baritz-expansion}),  namely 
$$
x^{-p} J_p(x) - \sum_{k=0}^N  \frac{ r^{-(p+k) } J_{p+k}(r) } { 2^k k! } (r^2 - x^2 )^k 
=  (r^2 - x^2 )^{N+1}  
    \sum_{k=N+1}^\infty  \frac{ r^{-(p+k) } J_{p+k}(r) } { 2^k k! } (r^2 - x^2 )^{k-N-1} .  
$$
Let  $ j_{n,1} $  denote the first positive zero of  $ J_n(x) $.  It is well known that
$ j_{n,1} < j_{n+1,1} $  and that  $ J_n(x) > 0 $  in  $ (0, j_{n,1} ) $.  See Figure 1.
\begin{center}
    \includegraphics[scale=0.6]{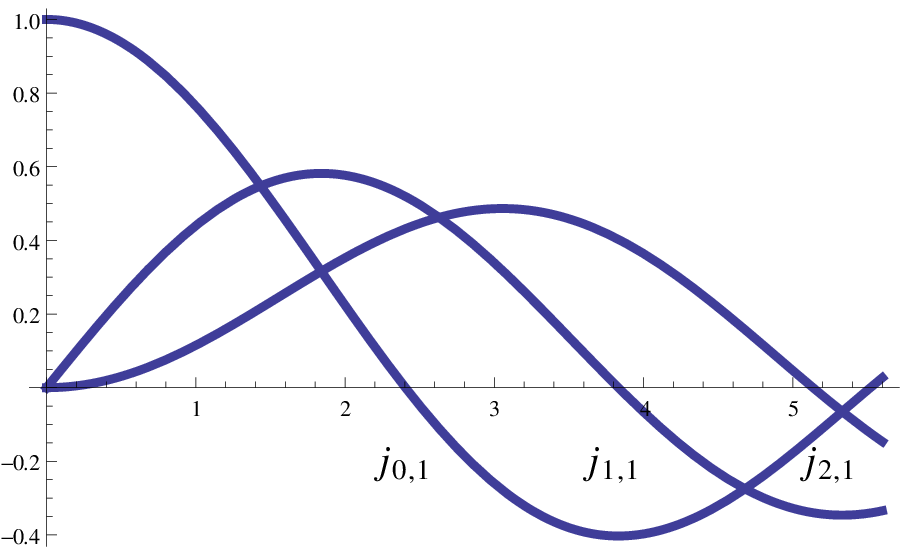}
\end{center}
\centerline{ Figure 1. The first zeros of $ J_0(x), \,  J_1(x)$ and $ J_2(x) $ }
If  $r$  satisfies  
	$  0 < r < j_{p+1,1} $, 
all the coefficients in the tail series  
$$ 
	R_N(x) = \sum_{k=N+1}^\infty  \frac{ r^{-(p+k) } J_{p+k}(r) } { 2^k k! } (r^2 - x^2 )^{k-N-1} 
$$
are positive for all  $ N $, hence 
	$  R_N(0)  \ge  R_N(x)  \ge  R_N(r) $. 
So, in the left hand side of  (\ref{eq:Baricz}), 
$ \alpha_p = R_N(r)=  \dfrac{ r^{-(p+N+1) } J_{p+N+1}(r) } { 2^{N+1} (N+1)! } $ ;
and on the right hand side of  (\ref{eq:Baricz}), 
$$  \beta_p = R_N(0)  
  	 = r^{ -2(N+1) } \sum_{k=N+1}^\infty  \frac{ r^{ k - p } J_{p+k}(r) } { 2^k k! }  \ . 
$$
$ \beta_p $  may be written as a finite sum.  
As  $ \lim_{x \to 0} x^{-p} J_p(x) = 1/ 2^{p} \Gamma(p+1) $,  we get from  
(\ref{Baritz-expansion})  for  $ x=0 $  that
$$    
   \frac{1}{2^p \Gamma(p+1) } = \sum_{k=0}^\infty \frac{ r^{ k - p } J_{p+k}(r) } { 2^k k! }  \ . 
$$
Thus,
$$  \beta_p =  \frac
   {  1 / 2^{p} \Gamma(p+1) - \sum_{k=0}^N   { r^{ k - p } J_{p+k}(r) } / { 2^k k! }  }
   { r^{ 2(N+1) }  }  \  .
$$

Inspired by a remark made by the referee, we suggest the following problem:  Find similar 
inequalities for other special functions, based on the above methods.



\setlength{\parindent}{0 mm}

Department of Mathematics, Technion --- I.I.T., Haifa 32000, Israel	\\ 
{\tt dova@tx.technion.ac.il}	\\ 
{\tt elias@tx.technion.ac.il}

\end{document}